# Some conjectures about q-Fibonacci polynomials

*Johann Cigler*

Fakultät für Mathematik
Universität Wien
A-1090 Wien, Nordbergstraße 15

Johann Cigler@univie.ac.at

## Abstract

In this paper we state some conjectures about $q-$Fibonacci polynomials which for $q=1$ reduce to well-known results about Fibonacci numbers and Fibonacci polynomials.

## 1. Introduction

The Fibonacci numbers $F_n$ are defined by the recurrence relation

$$F_n = F_{n-1} + F_{n-2} \tag{1.1}$$

with initial values $F_0 = 0$ and $F_1 = 1$.

The powers $F_n^k$, $k = 1, 2, 3, \cdots,$ satisfy the recurrence relation

$$\sum_{j=0}^{k+1}(-1)^{\binom{j+1}{2}} \left\langle \begin{matrix} k+1 \\ j \end{matrix} \right\rangle F_{n-j}^k = 0, \tag{1.2}$$

where $\left\langle \begin{matrix} n \\ k \end{matrix} \right\rangle = \dfrac{\prod_{i=0}^{k-1} F_{n-i}}{\prod_{i=1}^{k} F_i}$ is a so called fibonomial coefficient.

E.g. the squares of the Fibonacci numbers satisfy the recurrence $F_n^2 - 2F_{n-1}^2 - 2F_{n-2}^2 + F_{n-3}^2 = 0.$

The triangle of Fibonomial coefficients ( see A010048 or A055870 in the On-Line Encyclopedia of Integer Sequences [8] ) begins with

```
            1
          1   1
        1   1   1
      1   2   2   1
    1   3   6   3   1
  1   5  15  15   5   1
```

The Fibonacci polynomials $F_n(x,s)$ are defined by the recurrence relation

$$F_n(x,s) = xF_{n-1}(x,s) + sF_{n-2}(x,s) \tag{1.3}$$

with initial values $F_0(x,s) = 0$ and $F_1(x,s) = 1$.

The first terms of this sequence are $0, 1, x, x^2 + s, x^3 + 2sx, x^4 + 3sx^2 + s^2, \cdots$.

The powers $F_n^k(x,s)$, $k = 1, 2, 3, \cdots$, satisfy the recurrence relation

$$\sum_{j=0}^{k+1} (-1)^{\binom{j+1}{2}} s^{\binom{j}{2}} \left\langle{k+1 \atop j}\right\rangle(x,s) F_{n-j}^k(x,s) = 0, \tag{1.4}$$

where the polynomial fibonomial coefficients are defined by

$$\left\langle{n \atop k}\right\rangle(x,s) = \frac{\prod_{i=0}^{k-1} F_{n-i}(x,s)}{\prod_{i=1}^{k} F_i(x,s)}. \tag{1.5}$$

E.g. for $k = 2$ we get the recurrence relation
$F_n(x,s)^2 - (x^2 + s) F_{n-1}(x,s)^2 - s(x^2 + s) F_{n-2}(x,s)^2 + s^3 F_{n-3}(x,s)^2 = 0.$

The simplest proof of these facts depends on the Binet formula

$$F_n(x,s) = \frac{\alpha^n - \beta^n}{\alpha - \beta}, \tag{1.6}$$

where

$$\alpha = \frac{x + \sqrt{x^2 + 4s}}{2}, \beta = \frac{x - \sqrt{x^2 + 4s}}{2}. \tag{1.7}$$

From (1.6) it is clear that $F_n(x,s)^k$ is a linear combination of $\alpha^{(k-j)n} \beta^{jn}$, $0 \leq j \leq k$. Let $U$ be the shift operator $Uh(n) = h(n-1)$. The sequences $\left(\alpha^{(k-j)n} \beta^{jn}\right)_{n \geq 0}$ satisfy the recurrence relation $(1 - \alpha^{k-j} \beta^j U)(\alpha^{(k-j)n} \beta^{jn}) = 0$.

Since the operators $1 - \alpha^{k-j} \beta^j U$ commute we get

$$\left(\prod_{j=0}^{n} (1 - \alpha^{k-j} \beta^j U)\right) F_n(x,s)^k = 0. \tag{1.8}$$

As has been observed by L. Carlitz [2] we can now apply the $q$-binomial theorem (cf. e.g. [4])

$$\prod_{j=0}^{n-1} (1 - q^j x) = \sum_{k=0}^{n} (-1)^k q^{\binom{k}{2}} \begin{bmatrix} n \\ k \end{bmatrix} x^k. \tag{1.9}$$

Here $\begin{bmatrix} n \\ k \end{bmatrix} = \begin{bmatrix} n \\ k \end{bmatrix}(q) = \frac{(1-q^n)(1-q^{n-1}) \cdots (1-q^{n-k-1})}{(1-q)(1-q^2) \cdots (1-q^k)}$ denotes a $q$-binomial coefficient.

For $q = \dfrac{\beta}{\alpha}$ we get $\begin{bmatrix} n \\ k \end{bmatrix} = \alpha^{k^2 - nk} \left\langle \begin{matrix} n \\ k \end{matrix} \right\rangle (x,s)$.

This implies

$$\prod_{j=0}^{k}(1 - \alpha^{k-j}\beta^j U) = \prod_{j=0}^{k}\left(1 - \left(\frac{\beta}{\alpha}\right)^j (\alpha^k U)\right) = \sum_{j=0}^{k}(-1)^j \left(\frac{\beta}{\alpha}\right)^{\binom{j}{2}} \alpha^{j^2-(k+1)j} \left\langle \begin{matrix} k+1 \\ j \end{matrix} \right\rangle (x,s) \alpha^{kj} U^j$$

$$= \sum_{j=0}^{k}(-1)^j (\alpha\beta)^{\binom{j}{2}} \left\langle \begin{matrix} k+1 \\ j \end{matrix} \right\rangle (x,s) U^j = \sum_{j=0}^{k}(-1)^{\binom{j+1}{2}} s^{\binom{j}{2}} \left\langle \begin{matrix} k+1 \\ j \end{matrix} \right\rangle (x,s) U^j.$$

By applying this operator to $F_n(x,s)^k$ we get (1.4).

The Lucas polynomials
$$L_n(x,s) = \alpha^n + \beta^n \tag{1.10}$$

satisfy the same recurrence $L_n(x,s) = xL_{n-1}(x,s) + sL_{n-2}(x,s)$ with initial values $L_0(x,s) = 2$ and $L_1(x,s) = x$.

Therefore for $2j < k$ the coefficients of the product
$$(z - \alpha^{k-j}\beta^j)(z - \alpha^j \beta^{k-j}) = z^2 - (\alpha\beta)^j (\alpha^{k-2j} + \beta^{k-2j})z + (\alpha\beta)^k = z^2 - (-s)^j L_{k-2j}(x,s) z + (-s)^k$$
are polynomials in $x$ and $s$ with integer coefficients. The same is true if $k = 2j$. In this case we have $z - \alpha^{2j-j}\beta^j = z - (\alpha\beta)^j = z - (-s)^j$.

This implies that the coefficients of

$$\prod_{j=0}^{k}(z - \alpha^{k-j}\beta^j) = \sum_{j=0}^{k}(-1)^{\binom{j+1}{2}} s^{\binom{j}{2}} \left\langle \begin{matrix} k+1 \\ j \end{matrix} \right\rangle (x,s) z^{k+1-j} \tag{1.11}$$

are polynomials in $x$ and $s$ with integer coefficients.

The matrix of the polynomials $\left\langle \begin{matrix} k \\ j \end{matrix} \right\rangle (x,s)$ begins with

$$\begin{pmatrix} 1 & 0 & 0 & 0 & 0 & 0 \\ 1 & 1 & 0 & 0 & 0 & 0 \\ 1 & x & 1 & 0 & 0 & 0 \\ 1 & s+x^2 & s+x^2 & 1 & 0 & 0 \\ 1 & x(2s+x^2) & (s+x^2)(2s+x^2) & x(2s+x^2) & 1 & 0 \\ 1 & s^2+3sx^2+x^4 & (2s+x^2)(s^2+3sx^2+x^4) & (2s+x^2)(s^2+3sx^2+x^4) & s^2+3sx^2+x^4 & 1 \end{pmatrix}$$

**Remark**

The polynomial $\prod_{j=1}^{n}(z - \alpha^{j-1}\beta^{n-j}) = \sum_{j=0}^{n}(-1)^{\binom{j+1}{2}} s^{\binom{j}{2}} \left\langle \begin{matrix} n \\ j \end{matrix} \right\rangle (x,s) z^{n-j}$ also appears as the characteristic polynomial $\det(zI_n - a(n))$ of the matrix

$$a(n) = \left( \binom{i-1}{n-j} x^{i+j-n-1} s^{n-j} \right)_{i,j=1}^{n}.$$

For the special case $x = s = 1$ this conjecture by V.E. Hoggatt has been proved by L. Carlitz [2].

H. Prodinger [7] has given a simple proof that the eigenvectors $u(n, j)$ which belong to the eigenvalues $\lambda_j = \alpha^{j-1} \beta^{n-j}$, $1 \le j \le n$, are

$$u(n, j) = \begin{pmatrix} u(n,1,j) \\ u(n,2,j) \\ \vdots \\ u(n,n,j) \end{pmatrix}$$

with $u(n,i,j) = \sum_{k=1}^{j} (-s)^{i-k} \binom{i-1}{k-1} \binom{n-i}{j-k} \alpha^{2k-i-1}.$

## 2. Recurrence relations for powers of q-Fibonacci polynomials

The (Carlitz-) $q$–Fibonacci polynomials $f(n, x, s)$ are defined by

$$f(n, x, s) = xf(n-1, x, s) + q^{n-2} sf(n-2, x, s) \tag{2.1}$$

with initial values $f(0, x, s) = 0$, $f(1, x, s) = 1$ (cf. [3],[5]).
The first values are

0, 1, x, $q s + x^2$, $q s x + q^2 s x + x^3$, $q^4 s^2 + q s x^2 + q^2 s x^2 + q^3 s x^2 + x^4$.

An explicit expression is

$$f(n, x, s) = \sum_{k \le n-1} \begin{bmatrix} n-1-k \\ k \end{bmatrix} q^{k^2} x^{n-1-2k} s^k. \tag{2.2}$$

If we change $q \to \dfrac{1}{q}$ and then $s \to q^{n-1} s$ we get

$$f(n-k, x, q^{-k} s) \to f(n-k, x, q^{1-n} s) \to f(n, x, s).$$

Therefore each identity

$$g(x, s, q, f(n, x, s), f(n-1, x, s), f(n-2, x, s), \cdots) = 0 \tag{2.3}$$

is equivalent with the identity

$$g(x, q^{n-1} s, q^{-1}, f(n, x, s), f(n-1, x, qs), f(n-2, x, q^2 s), \cdots) = 0. \tag{2.4}$$

As a special case we get the well-known fact that (2.1) is equivalent with

$$f(n, x, s) = xf(n-1, x, qs) + qsf(n-2, x, q^2 s). \tag{2.5}$$

The definition of the $q-$Fibonacci polynomials can be extended to all integers such that the recurrence (2.1) remains true. We then get (cf. [5])

$$f(-n, x, s) = (-1)^{n-1} q^{\binom{n+1}{2}} \frac{f(n, x, q^{-n} s)}{s^n}. \tag{2.6}$$

We are now looking for a $q-$analog of (1.4).
Computer experiments suggest the following

**Conjecture 1**

*Let*

$$\left\langle {k \atop j} \right\rangle (x, s, q) = \frac{\prod_{i=1}^{k} f(i, x, s)}{\prod_{i=1}^{j} f(i, x, q^{j-i} s) \prod_{i=1}^{k-j} f(i, x, q^j s)}. \tag{2.7}$$

*Then the following recurrence relation holds for all $n \in \mathbb{Z}$:*

$$\sum_{j=0}^{k+1} (-1)^{\binom{j+1}{2}} s^{\binom{j}{2}} q^{\frac{j(j-1)(2j-1)}{6}} \left\langle {k+1 \atop j} \right\rangle (x, s, q) f(n-j, x, q^j s)^k = 0. \tag{2.8}$$

*This relation is equivalent to*

$$\sum_{j=0}^{k+1} (-1)^{\binom{j+1}{2}} s^{\binom{j}{2}} q^{(n-1)\binom{j}{2} - \frac{j(j-1)(2j-1)}{6}} fibo(k+1, x, s) f(n-j, x, s)^k = 0 \tag{2.9}$$

*with*

$$fibo(k+1, x, s) = \left\langle {k+1 \atop j} \right\rangle (x, q^{n-1} s, q^{-1}) = \frac{\prod_{i=1}^{k} f(i, x, q^{n-i} s)}{\prod_{i=1}^{j} f(i, x, q^{n-j} s) \prod_{i=1}^{k-j} f(i, x, q^{n-i-j} s)}. \tag{2.10}$$

This conjecture is trivially true for $k = 1$.
It should be noted that $f(n, x, s)$ is a $q-$holonomic sequence of polynomials.
A sequence $(a(n))$ is $q-$holonomic if there exist nontrivial polynomials $p_0, \cdots, p_r$ such that
$p_0(q^n) a(n) + p_1(q^n) a(n-1) + \cdots + p_r(q^n) a(n-r) = 0$.
If $(a(n))$ is $q-$holonomic then for each $k \in \mathbb{N}$ the sequence $\left( a(n)^k \right)$ also is $q-$holonomic.
Therefore it is clear that $\left( f(n, x, s)^k \right)$ has a $q-$holonomic recurrence. With the help of the Mathematica package qGeneratingFunction by Christoph Koutschan (RISC Linz) I have computed recursion (2.9) for small values of $k$.

For $k = 2$ (2.9) reduces to

$$f(n,x,s)^2 - (x^2 + q^{n-2}s)f(n-1,x,s)^2 - q^{n-2}s(x^2 + q^{n-2}s)f(n-2,x,s)^2 + q^{3n-8}s^3 f(n-3,x,s)^2 = 0$$

and (2.8) to

$$f(n,x,s)^2 - (x^2 + qs)f(n-1,x,qs)^2 - qs(x^2 + qs)f(n-2,x,q^2s)^2 \\ + q^5 s^3 f(n-3,x,q^3s)^2 = 0. \quad (2.11)$$

To prove this we use the $q$–Euler-Cassini formula (cf. [5], Cor. 2.2, see also (3.7))

$$f(k-1,x,qs)f(n+k,x,s) - f(k,x,s)f(n+k-1,x,qs) = (-1)^k q^{\binom{k}{2}} s^{k-1} f(n,x,q^k s). \quad (2.12)$$

For each $k \in \mathbb{N}$ it gives a representation of $f(n-k, x, q^k s)$ as a linear combination of $f(n,x,s)$ and $f(n-1,x,qs)$:

$$f(n-k, x, q^k s) = \frac{1}{v(k)}\left(f(k-1,x,qs)f(n,x,s) - f(k,x,s)f(n-1,x,qs)\right) \quad (2.13)$$

with

$$v(k) = (-1)^k q^{\binom{k}{2}} s^{k-1}.$$

If we set $f(n,x,s) = a$, $f(n-1,x,qs) = b$ we get

$$f(n-2, x, q^2 s) = \frac{1}{qs}(a - xb) \text{ and } f(n-3, x, q^3 s) = \frac{-1}{q^3 s^2}\left(xa - (x^2 + qs)b\right).$$

Therefore (2.11) becomes

$$a^2 - (x^2 + qs)b^2 - \frac{qs(x^2 + qs)}{(qs)^2}(a - xb)^2 + \frac{q^5 s^3}{(-q^3 s^2)^2}\left(xa - (x^2 + qs)b\right)^2 = 0.$$

This is evident since all the coefficients of $a^2, ab, b^2$ vanish.

In the same way I have verified the conjecture for other small values of $k$.

### 3. Recurrence relations for subsequences
In the classical case $q = 1$ formula (1.4) can be generalized to

$$\sum_{j=0}^{k+1}(-1)^{j+\ell\binom{j}{2}} s^{\ell\binom{j}{2}} \left\langle \begin{matrix} k+1 \\ j \end{matrix} \right\rangle(\ell, x, s)\left(F_{\ell(n-j)}(x,s)\right)^k = 0 \quad (3.1)$$

with

$$\left\langle \begin{matrix} k \\ j \end{matrix} \right\rangle(\ell, x, s) = \frac{\prod_{i=0}^{j-1} F_{(k-i)\ell}(x,s)}{\prod_{i=1}^{j} F_{i\ell}(x,s)}. \quad (3.2)$$

This follows from the same argument as above by changing $\alpha \to \alpha^\ell, \beta \to \beta^\ell$.

An analogous formula seems to hold in the general case.

## Conjecture 2

*For all integers $k \geq 1$ and $\ell \geq 1$ the following recurrence holds for all $n \in \mathbb{Z}$:*

$$\sum_{j=0}^{k+1}(-1)^{j+\ell\binom{j}{2}}\left(q^{\frac{(4j+1)\ell-3}{6}}s\right)^{\ell\binom{j}{2}}\left\langle {k+1 \atop j} \right\rangle(\ell,x,s,q)f(\ell(n-j),x,q^{\ell j}s)^k = 0 \qquad (3.3)$$

with

$$\left\langle {m \atop j} \right\rangle(\ell,x,s,q) = \frac{\prod_{i=1}^{m}f(\ell i,x,s)}{\prod_{i=1}^{j}f(\ell i,x,q^{\ell(j-i)}s)\prod_{i=1}^{m-j}f(\ell i,x,q^{\ell j}s)}. \qquad (3.4)$$

Let us first consider the case $k=1$.
Here (3.3) reduces to

$$f(\ell n,x,s) - \frac{f(2\ell,x,s)}{f(\ell,x,q^\ell s)}f(\ell(n-1),x,q^\ell s) + (-1)^\ell q^{\frac{\ell(3\ell-1)}{2}}s^\ell\frac{f(\ell,x,s)}{f(\ell,x,q^\ell s)}f(\ell(n-2),x,q^{2\ell}s) = 0. \qquad (3.5)$$

For $q=1$ this is the well-known recurrence

$$F_{\ell n}(x,s) - L_\ell(x,s)F_{\ell(n-1)}(x,s) + (-s)^\ell F_{\ell(n-2)}(x,s) = 0. \qquad (3.6)$$

For in this case we have $\dfrac{F_{2\ell}(x,s)}{F_\ell(x,s)} = \dfrac{\alpha^{2\ell}-\beta^{2\ell}}{\alpha^\ell-\beta^\ell} = \alpha^\ell + \beta^\ell = L_\ell(x,s)$.

In order to prove (3.5) we prove the slightly more general formula

$$\det\begin{pmatrix} f(N+(m+1)\ell,x,s) & f((m+1)\ell,x,s) \\ f(N+m\ell,x,q^\ell s) & f(m\ell,x,q^\ell s) \end{pmatrix} = (-1)^{m\ell-1}s^{m\ell}q^{\frac{m\ell((m+2)\ell-1)}{2}}f(\ell,x,s)f(N,x,q^{(m+1)\ell}s). \qquad (3.7)$$

From $f(N+m\ell,x,q^\ell s) = xf(N+m\ell-1,x,q^\ell s) + q^{N-2+(m+1)\ell}sf(N+m\ell-2,x,q^\ell s)$
and $f(N+(m+1)\ell,x,s) = xf(N+(m+1)\ell-1,x,s) + q^{N-2+(m+1)\ell}sf(N+(m+1)\ell-2,x,s)$
we see that

$$g(N):=\det\begin{pmatrix} f(N+(m+1)\ell,x,s) & f((m+1)\ell,x,s) \\ f(N+m\ell,x,q^\ell s) & f(m\ell,x,q^\ell s) \end{pmatrix}$$

satisfies $g(N) = xg(N-1) + q^{N-2+(m+1)\ell}sg(N-2)$ and $g(0)=0$. Therefore
$g(N) = cf(N,x,q^{(m+1)\ell}s)$ for some constant $c$. To compute $c$ we set $N=-m\ell$. This gives
$g(-m\ell) = f(\ell,x,s)f(m\ell,x,q^\ell s) = cf(-m\ell,x,q^{(m+1)\ell}s)$ or
$c = -(-s)^{m\ell}q^{\frac{m\ell((m+2)\ell-1)}{2}}f(\ell,x,s)$. Therefore we get (3.7).

## Remark

For $k=1, m \to k-1$ (3.7) gives the $q$-Euler-Cassini formula (2.12).

If we let $n \to \infty$ in (2.2) with $x = 1$ we get

$$F(s) = \sum_{k \geq 0} \frac{q^{k^2}}{(1-q)(1-q^2)\cdots(1-q^k)} s^k. \tag{3.8}$$

If we let $n \to \infty$ in the $q$-Euler-Cassini formula we get

$$f(k-1,1,qs)F(s) - f(k,1,s)F(qs) = (-1)^k q^{\binom{k}{2}} s^{k-1} F(q^k s).$$

The special case $s = 1$ is closely connected with the famous Rogers-Ramanujan identities (cf. [1]).

For $k = 2$ formula (3.3) reduces to

$$f(\ell n, x, s)^2 - \frac{f(3\ell, x, s) f(2\ell, x, s)}{f(\ell, x, q^\ell s) f(2\ell, x, q^\ell s)} f(\ell n - \ell, x, q^\ell s)^2$$

$$+ (-1)^\ell q^{\frac{\ell(3\ell-1)}{2}} s^\ell \frac{f(3\ell, x, s) f(\ell, x, s)}{f(\ell, x, q^\ell s) f(\ell, x, q^{2\ell} s)} f(\ell n - 2\ell, x, q^{2\ell} s)^2 \tag{3.9}$$

$$+ (-1)^{\ell-1} s^{3\ell} q^{\frac{\ell(13\ell-3)}{2}} \frac{f(\ell, x, s) f(2\ell, x, s)}{f(\ell, x, q^{2\ell} s) f(2\ell, x, q^\ell s)} f(\ell n - 3\ell, x, q^{3\ell} s)^2.$$

This can be proved in the same way as in the special case $\ell = 1$ by using (3.5) in order to reduce all expressions $f(\ell n - j\ell, x, q^{j\ell} s)^2$ to $a = f(\ell n - 3\ell, x, q^{3\ell} s)$ and $b = f(\ell n - 2\ell, x, q^{2\ell} s)$.

For other small values of $k$ formula (3.3) can be verified in the same way. But I did not find a proof for all $k$.

## 4. A generalization of Cassini's identity

A well-known formula for Fibonacci numbers is Cassini's identity

$$\det \begin{pmatrix} F_n & F_{n-1} \\ F_{n-1} & F_n \end{pmatrix} = (-1)^{n-1}. \tag{4.1}$$

This can be generalized in the following way:

Let

$$Fac(n, s, \ell) = F_\ell(x, s) F_{2\ell}(x, s) \cdots F_{n\ell}(x, s). \tag{4.2}$$

Then

$$\det \left( F_{n+i-j}(x, s)^k \right)_{i,j=0}^k = (-1)^{\binom{k+1}{2}(n-k)} \prod_{\ell=0}^k \binom{k}{\ell} s^{\binom{k+1}{2}(n-k)+2\binom{k+1}{3}} \prod_{j=0}^{k-1} Fac(k-j, s, 1)^2. \tag{4.3}$$

As a special case we get

$$\det\begin{pmatrix} F_n^2 & F_{n-1}^2 & F_{n-2}^2 \\ F_{n+1}^2 & F_n^2 & F_{n-1}^2 \\ F_{n+2}^2 & F_{n+1}^2 & F_n^2 \end{pmatrix} = 2(-1)^n.$$

To prove this we use again Binet's formula:

$$\det\left(F_{n+i-j}(x,s)^k\right)_{i,j=0}^k = \det\left(\frac{\left(\alpha^{n+i-j} - \beta^{n+i-j}\right)^k}{(\alpha - \beta)^k}\right) = \left(\frac{1}{(\alpha-\beta)^{k^2+k}}\right)\det\left(\sum_{\ell=0}^k (-1)^\ell \binom{k}{\ell} \alpha^{(n+i-j)(k-\ell)} \beta^{(n+i-j)\ell}\right).$$

Let now $z(j,\ell) = \begin{pmatrix} \alpha^{(n-j)(k-\ell)} \beta^{(n-j)\ell} \\ \alpha^{(n+1-j)(k-\ell)} \beta^{(n+1-j)\ell} \\ \vdots \\ \alpha^{(n+k-j)(k-\ell)} \beta^{(n+k-j)\ell} \end{pmatrix}.$

Then we get

$$\det\left(\sum_{\ell=0}^k (-1)^\ell \binom{k}{\ell} \alpha^{(n+i-j)(k-\ell)} \beta^{(n+i-j)\ell}\right) = (-1)^{\lfloor \frac{k+1}{2} \rfloor} \prod_{\ell=0}^k \binom{k}{\ell} \sum_\pi \det\left(z(0,\pi(0)), z(1,\pi(1)), \cdots z(k,\pi(k))\right) \mathrm{sgn}(\pi),$$

where the sum is over all permutations $\pi$ of $\{0,\cdots,k\}$.

This implies

$$\det\left(F_{n+i-j}(x,s)^k\right)_{i,j=0}^k = \left(\frac{1}{(\alpha-\beta)^{k(k+1)}}\right)(-1)^{\lfloor \frac{k+1}{2} \rfloor} \prod_{\ell=0}^k \binom{k}{\ell} \left(\sum_\pi \mathrm{sgn}(\pi) \prod_{i=0}^k \alpha^{(n-i)(k-\pi(i))} \beta^{(n-i)\pi(i)}\right) d$$

with $d = \det\left(\left(\alpha^{k-j}\beta^j\right)^i\right).$

By the formula for Vandermonde's determinant we get

$$d = \det\left(\left(\alpha^{k-j}\beta^j\right)^i\right) = (-1)^{\binom{k+1}{2}} \left(\alpha^k - \alpha^{k-1}\beta\right)\left(\alpha^k - \alpha^{k-2}\beta^2\right)\cdots\left(\alpha^k - \beta^k\right)\left(\alpha^{k-1}\beta - \alpha^{k-2}\beta^2\right)$$
$$\cdots\left(\alpha^{k-1}\beta - \beta^k\right)\cdots\left(\alpha\beta^{k-1} - \beta^k\right)$$
$$= (-1)^{\binom{k+1}{2}} (\alpha\beta)^{\binom{k}{2}+\binom{k-1}{2}+\cdots+\binom{k}{2}} (\alpha - \beta)(\alpha^2 - \beta^2)\cdots(\alpha^k - \beta^k)(\alpha - \beta)\cdots(\alpha^{k-1} - \beta^{k-1})\cdots(\alpha - \beta)$$
$$= (-1)^{\binom{k+1}{2}} (-s)^{\binom{k+1}{3}} (\alpha-\beta)^{\binom{k+1}{2}} \prod_{j=0}^{k-1} Fac(k-j,s,1).$$

The sum $\sum_\pi \mathrm{sgn}(\pi) \prod_{i=0}^k \alpha^{(n-i)(k-\pi(i))} \beta^{(n-i)\pi(i)}$ equals $\det(c_{i,j})$ with $c_{i,j} = \left(\alpha^{k-j}\beta^j\right)^{n-i}$, which gives

$$\det(c_{i,j}) = d \prod_{j=0}^k \left(\alpha^{k-j}\beta^j\right)^{n-k}.$$

Therefore we finally get

$$\det\left(F_{n+i-j}(x,s)^k\right)_{i,j=0}^k = \left(\frac{1}{(\alpha-\beta)^{k(k+1)}}\right)(-1)^{\lfloor\frac{k+1}{2}\rfloor+\binom{k+1}{2}}\prod_{\ell=0}^k\binom{k}{\ell}d^2\prod_{j=0}^k\left(\alpha^{k-j}\beta^j\right)^{n-k}$$

$$= \left(\frac{1}{(\alpha-\beta)^{k(k+1)}}\right)\prod_{\ell=0}^k\binom{k}{\ell}d^2(-s)^{\binom{k+1}{2}(n-k)}d^2$$

$$= \left(\frac{1}{(\alpha-\beta)^{k(k+1)}}\right)(-1)^{\binom{k+1}{2}(n-k)}\prod_{\ell=0}^k\binom{k}{\ell}s^{\binom{k+1}{2}(n-k)+2\binom{k+1}{3}}(\alpha-\beta)^{2\binom{k+1}{2}}\prod_{j=0}^{k-1}Fac(k-j,s,1)^2.$$

A $q$-analog of Cassini's identity is

$$d(n,s) = \det\begin{pmatrix} f(n,x,s) & f(n-1,x,qs) \\ f(n+1,x,s) & f(n,x,qs) \end{pmatrix}$$

$$= \det\begin{pmatrix} xf(n-1,x,qs)+qsf(n-2,x,q^2s) & f(n-1,x,qs) \\ xf(n,x,qs)+qsf(n-1,x,q^2s) & f(n,x,qs) \end{pmatrix}$$

$$= \det\begin{pmatrix} qsf(n-2,x,q^2s) & f(n-1,x,qs) \\ qsf(n-1,x,q^2s) & f(n,x,qs) \end{pmatrix} = -qsd(n-1,qs)$$

which gives $d(n,s) = (-1)^{n-1}q^{\binom{n}{2}}s^{n-1}$.

This identity can be generalized to

$$\det\begin{pmatrix} f(n,x,s)^2 & f(n-1,x,qs)^2 & f(n-2,x,q^2s)^2 \\ f(n+1,x,s)^2 & f(n,x,qs)^2 & f(n-1,x,q^2s)^2 \\ f(n+2,x,s)^2 & f(n+1,x,qs)^2 & f(n,x,q^2s)^2 \end{pmatrix} = 2(-1)^n x^2 s^{3n-4} q^{\frac{(n+1)(3n-4)}{2}}. \qquad (4.4)$$

For $n=2$ this can be verified by computation. In the general case we have using (2.11)

$$d(n,s) = \det\begin{pmatrix} f(n,x,s)^2 & f(n-1,x,qs)^2 & f(n-2,x,q^2s)^2 \\ f(n+1,x,s)^2 & f(n,x,qs)^2 & f(n-1,x,q^2s)^2 \\ f(n+2,x,s)^2 & f(n+1,x,qs)^2 & f(n,x,q^2s)^2 \end{pmatrix}$$

$$= \det\begin{pmatrix} q^5s^3f(n-3,x,q^3s)^2 & f(n-1,x,qs)^2 & f(n-2,x,q^2s)^2 \\ q^5s^3f(n-2,x,q^3s)^2 & f(n,x,qs)^2 & f(n-1,x,q^2s)^2 \\ q^5s^3f(n-1,x,q^3s)^2 & f(n+1,x,qs)^2 & f(n,x,q^2s)^2 \end{pmatrix} = -q^5s^3d(n-1,qs).$$

This implies (4.4).

The same method applies for each small $n$. These results lead to the following

## Conjecture 3

Let $fac(n, s, \ell) = f(\ell, x, s) f(2\ell, x, s) \cdots f(n\ell, x, s)$.
Then for $n \in \mathbb{Z}$ and $k \geq 1$ the following identity holds:

$$\det\left(f(n+i-j, q^j s)^k\right)_{i,j=0}^{k}$$
$$= \left((-1)^{(n-k)\binom{k+1}{2}} \prod_{j=0}^{k}\binom{k}{j}\right)\left(q^{\frac{n+k-1}{2}} s\right)^{2\binom{k+1}{3}+\binom{k+1}{2}(n-k)} \prod_{j=0}^{k-1} fac(k-j, q^j s, 1) \prod_{j=0}^{k-1} fac(k-j, q^{n+j} s, 1). \quad (4.5)$$

For $q = 1$ the same proof as above gives

$$\det\left(F_{\ell(n+i-j)}(x, s)^k\right)_{i,j=0}^{k} = (-1)^{\ell\binom{k+1}{2}(n-k)} \left(\prod_{m=0}^{k}\binom{k}{m}\right) s^{\ell\left(\binom{k+1}{2}(n-k)+2\binom{k+1}{3}\right)} \prod_{j=0}^{k-1} Fac((k-j), s, \ell)^2.$$

This leads to

## Conjecture 4

For $n \in \mathbb{Z}$ and $\ell, k \geq 1$ the following identity holds:

$$\det\left(f(\ell(n+i-j), x, q^{\ell j} s)^k\right)_{i,j=0}^{k} = (-1)^{\ell\binom{k+1}{2}(n-k)} \left(\prod_{m=0}^{k}\binom{k}{m}\right) \left(q^{\frac{\ell(n+k)-1}{2}} s\right)^{\ell\left(\binom{k+1}{2}(n-k)+2\binom{k+1}{3}\right)}$$
$$\prod_{j=0}^{k-1} fac(k-j, q^{\ell j} s, \ell) \prod_{j=0}^{k-1} fac(k-j, q^{\ell(n+j)} s, \ell). \quad (4.6)$$

The simplest special cases are

$$\det\begin{pmatrix} f(n\ell, x, s) & f((n-1)\ell, x, q^\ell s) \\ f((n+1)\ell, x, s) & f(n\ell, x, q^\ell s) \end{pmatrix} = (-1)^{(n-1)\ell} \left(q^{\frac{\ell(n+1)-1}{2}} s\right)^{(n-1)\ell} f(\ell, x, s) f(\ell, x, q^{n\ell} s) \quad (4.7)$$

and

$$\det\begin{pmatrix} f(\ell n, x, s)^2 & f(\ell(n-1), x, q^\ell s)^2 & f(\ell(n-2), x, q^{2\ell} s)^2 \\ f(\ell(n+1), x, s)^2 & f(\ell n, x, q^\ell s)^2 & f(\ell(n-1), x, q^{2\ell} s)^2 \\ f(\ell(n+2), x, s)^2 & f(\ell(n+1), x, q^\ell s)^2 & f(\ell n, x, q^{2\ell} s)^2 \end{pmatrix}$$
$$= 2(-1)^{n\ell} \left(q^{\frac{\ell(n+2)-1}{2}} s\right)^{\ell(3n-4)} f(2\ell, x, s) f(2\ell, x, q^{n\ell} s) f(\ell, x, s) f(\ell, x, q^{\ell} s) f(\ell, x, q^{n\ell} s) f(\ell, x, q^{(n+1)\ell} s). \quad (4.8)$$

Formula (4.7) is a special case of (3.7). If we let $N = -\ell$ in (3.7) we get

$$\det\begin{pmatrix} f(m\ell, x, s) & f((m+1)\ell, x, s) \\ f((m-1)\ell, x, q^\ell s) & f(m\ell, x, q^\ell s) \end{pmatrix} = (-1)^{m\ell-1} s^{m\ell} q^{\frac{m\ell((m+2)\ell-1)}{2}} f(\ell, x, s) f(-\ell, x, q^{(m+1)\ell} s),$$

which in turn implies (4.7).

A proof of (4.8) may be found in [6].

In order to prove (4.6) it suffices to prove it for $n = k$ and then use induction by applying (3.3) to the first column.
So if (3.3) is already known it would suffice to prove

$$\det\left(f(\ell(k+i-j),x,q^{\ell j}s)^k\right)_{i,j=0}^{k} = (-1)^{\ell\binom{k+1}{2}(n-k)} \left(\prod_{j=0}^{k}\binom{k}{j}\right) q^{\ell(2\ell k-1)\binom{k+1}{3}} s^{2\ell\binom{k+1}{3}}$$

$$\prod_{j=0}^{k-1} fac(k-j,q^{\ell j}s,\ell) \prod_{j=0}^{k-1} fac(k-j,q^{\ell(k+j)}s,\ell).$$

## Remarks

The conjectures which I have stated in this paper are very simple and natural. Therefore I am sure that they are correct. Although there is no $q-$analog of the Binet formula, the conjectured formulae are so similar to their counterparts in the classical case (which heavily depend on Binet's formula), that there must be some reason for this fact. Perhaps this can lead to a proof. Especially interesting in my opinion is the occurrence of the factor $\prod_{j=0}^{k}\binom{k}{j}$ in the determinant evaluations. In fact I have computed the determinants $\det\left(f(k+i-j,q^j s)^k\right)_{i,j=0}^{k}$
for $1 \leq k \leq 5$.
```
{1, 2 q³ s² x², 9 q²⁰ s⁸ x⁴ (q s + x²) (q⁴ s + x²) , 96 q⁷⁰ s²⁰ x⁸ (q s + x²)
    (q² s + x²) (q s + q² s + x²) (q⁵ s + x²) (q⁶ s + x²) (q⁵ s + q⁶ s + x²) ,
  2500 q¹⁸⁰ s⁴⁰ x¹² (q s + x²) (q² s + x²) (q s + q² s + x²) (q³ s + x²) (q² s + q³ s + x²)
    (q⁶ s + x²) (q⁷ s + x²) (q⁶ s + q⁷ s + x²) (q⁸ s + x²) (q⁷ s + q⁸ s + x²)
    (q⁴ s² + q s x² + q² s x² + q³ s x² + x⁴) (q¹⁴ s² + q⁶ s x² + q⁷ s x² + q⁸ s x² + x⁴) }
```
I wondered how the sequence $\{1,2,9,96,2500,\cdots\}$ can be continued. So I asked Sloane's On-Line Encyclopedia of Integer Sequences [8] and found it there as sequence A001142. Only then did I compute the corresponding determinants for $q = 1$, where these numbers occur in a natural way. By the way I guess that all mentioned results for $q = 1$ are well known, but I don't know any references. So I would be very grateful to get some historical remarks or hints to the literature.